\documentclass[review]{elsarticle}

\usepackage{lineno,hyperref}
\usepackage{amsfonts}
\usepackage{amsmath}
\usepackage{amsxtra}
\usepackage{stmaryrd}
\usepackage{mathabx}
\usepackage{pdflscape}

\journal{Journal of \LaTeX\ Templates}

    \newtheorem{theorem}{Theorem}[section]
    \newtheorem{lemma}[theorem]{Lemma}
    
    \newtheorem{corollary}[theorem]{Corollary}

    \newenvironment{proof}[1][Proof]{\begin{trivlist}
    \item[\hskip \labelsep {\bfseries #1}]}{\end{trivlist}}
    \newenvironment{definition}[1][Definition]{\begin{trivlist}
    \item[\hskip \labelsep {\bfseries #1}]}{\end{trivlist}}
    \newenvironment{example}[1][Example]{\begin{trivlist}
    \item[\hskip \labelsep {\bfseries #1}]}{\end{trivlist}}
    \newenvironment{remark}[1][Remark]{\begin{trivlist}
    \item[\hskip \labelsep {\bfseries #1}]}{\end{trivlist}}

\begin{document}

\begin{frontmatter}

\title{On the number of antichains of sets in a finite universe}

\author{Patrick De Causmaecker\fnref{myfootnote} Stefan De Wannemacker}
\address{KU Leuven, Department of Computer Science, KULAK, CODeS \& iMinds-ITEC}


\ead{patrick.decausmaecker@kuleuven-kulak.be}


\begin{abstract}
Properties of intervals in the lattice of antichains of subsets of a  universe of finite size are investigated.
New objects and quantities in this lattice are defined.
Expressions and numerical values are deduced for the number of connected antichains and the number of fully distinguishing antichains.
The latter establish a connection with Stirling numbers of the second kind.
Decomposition properties of intervals in the lattice of antichains are proven.
A new operator allowing partitioning the full lattice in intervals derived from lower dimensional sub-lattices is introduced.
Special posets underlying an interval of antichains are defined.
The poset allows the derivation of a powerful formula for the size of an interval.
This formula allows computing intervals in the six dimensional space.
Combinatorial coefficients allowing another decomposition of the full lattice are defined.
In some specific cases, related to connected components in graphs, these coefficients can be efficiently computed.
This formula allows computing the size of the lattice of order 8 efficiently.
This size is the number of Dedekind of order 8, the largest one known so far.
\end{abstract}


\end{frontmatter}


\section{Introduction}
\label{sec:introduction}
The Dedekind number counts the number of antichains of subsets of an $n$-element set or the number of elements in a free distributive lattice on $n$  generators. Equivalently it counts number of monotonic functions on the subsets of a finite set of $n$ elements \cite{DEDEKIND,SLOANE}.
In 1969, Kleitman \cite{KLEITMAN} obtained an upper bound on the logarithm of the Dedekind number which was improved by Kleitman and Markowsky \cite{KLEITMAN_MARKOWSKY} in 1975
\[ log |M(n)| < (1+O((log n)/n)) \binom{n}{\lfloor n/2 \rfloor}.\]
In 1981, Korshunov \cite{KORSHUNOV} used a sophisticated approach to give assymptotics for the Dedekind number itself. All these proofs were simplified by Kahn \cite{KAHN} in 2002 using an 'Entropy' approach.
Finding a closed-form expression is a very hard problem, also known as Dedekind's problem and exact values have been found only for $n \leq 8$ \cite{CHURCH1940,CHURCH1965,WARD, WIEDEMANN}. This is sequence A000372 in Sloane's Online Encyclopedia of Integer Sequences \cite{SLOANE}. Recent attempts to improve the computational time are e.g. \cite{BAKOEV, FYDITEK}.
Recent related work on inequivalent antichains is in \cite{TAMON2014}.
The Dedekind numbers count the number of monotone boolean function as well as the number of antichains or the number of Sperner families.
The methods in the present paper are most naturally and easily understood in terms of antichains of subsets of a finite set.
Section \ref{sec:notationsandconventions} sets notations and conventions accordingly.

The remaining sections of the paper present several new expressions and recursion formulae for Dedekind numbers.
All formulae are based on partitions of the lattice of antichains of subsets of a finite set (referred to as 'the universe').
The antichains can be ordered in a natural way and lattice operators meet and join can be defined.
An important role in this paper is played by intervals in this lattice and the posets that underly them.
In Section \ref{sec:intervalposets}, four observations are stated in the form of lemmata to be used in the proofs of the theorems in the following sections.
The observations are a characterisation of posets underlying an interval (Lemma \ref{lem:posetcharacteristic}), 
a decomposition lemma for such posets (Lemma \ref{lem:posetdecomposition})
and two lemmata establishing isomorphism of posets underlying intervals and allowing dimensionality reduction (Lemmata \ref{lem:posetredundancy} and 
\ref{lem:posetindistinguishability}).
Section \ref{sec:intervals} lists three additional observations on intervals.
Lemma \ref{lem:intervalboundaryrendundancy} shows how to remove common elements from the top and  bottom antichains of an interval.
Lemma \ref{lem:intervalredundancy} allows to reduce top and bottom even further in case specific structures are recognised.
Finally, Lemma \ref{lem:graphdecomposition} associates a graph with the top antichain of an interval to decompose it as the direct join of intervals associated with the connected components of the graph.
Sections \ref{sec:partitioningbyintervals} and \ref{sec:decomposebyproperty} build on all these observations to prove decomposition and recursion formulae.
Two classes of formulae are distinguished: those based on partitions in intervals (Section \ref{sec:partitioningbyintervals}) and those based on properties of antichains (Section \ref{sec:decomposebyproperty}).
A first partition in intervals is related to any antichain in the lattice and uses the concept of the largest non dominating antichain (Theorem \ref{the:largestnondominatingdecomposition}).
A second class of partitions in intervals relies on a new operator. It allows a Cartesian-like decomposition of the lattice where the role of coordinates is taken by antichains in a lower dimension (Theorem \ref{the:intervalintervaldecomposition}).
In Section \ref{sec:decomposebyproperty}, a first property of antichains is the size of the sets in the antichain.
By constraining the antichain in an interval on subchains containing sets of specific sizes, it is possible to derive compact formulae for the number of antichains satisfying these constraints. In case the sizes differ by two only, these numbers are powers of two. Summing over all constraints gives the total number (Theorem \ref{the:intervallatticesize}).
The number of partitions of the universe an antichain can distinguish is the second property that can be used. A set is distinguished by an antichain if for each element outside this set, there is at least one set in the antichain containing the set and not this element. The smallest distinguished sets in an antichain form a partition of its span and Subsection \ref{sub:dedekindstirling2} uses the fact that the structure of such an antichain does not rely on the size of the smallest distinguished sets but solely on their number. The result is an expression for the number of antichains relying on Stirling numbers of the second kind (Theorem \ref{the:dedekindstirling2}).
Connectedness in antichains is the property studied in Section \ref{sub:connectedness}.
Two elements of the universe are connected in an antichain if there is a path from one element to another through the sets of the antichain.
Again, this property defines partitions of the universe and an expression for the total number of antichains could be derived (Theorem \ref{the:connectedness}).
Finally, in Section \ref{sub:pcoefficients}, an expansion of Dedekind numbers is given in terms of products of sizes of intervals in a lower dimensional space. 
The coefficients in the expansion are the number of solutions to a specific system of equations in the lattice (Theorem \ref{the:pcoefficientstheorem}).
They take a particularly easy form, related to connected subgraphs of a graph, when the difference in dimension is two.

All expressions in this paper have been implemented in java, an executable archive file is available at \cite{CODESREPORTS}.
Some of the implementations, such as the ones for Theorem \ref{the:pcoefficientstheorem}, use permutation symmetry and multithreading to speed up the computation. Computing the Dedekind number of order eight on a 3 Ghz computer with four cores takes about 8 hours.

Part of this work and early versions of some of the theorems were presented at CLA 2013 \cite{PDCSDW2013}.
This paper extends on these results and relates the theorems to the structure of the underlying posets introduced in Section \ref{sec:intervalposets}.

\section{Notations and conventions}
\label{sec:notationsandconventions}
Sets of subsets of a finite set are called antichains if no two of those subsets are in a subset-superset relationship.
If $N \subset \mathbb{N}$ is a finite set of numbers with $|N| = n$, the set of antichains in $N$ is
$
\mathcal{A}_N = \{ \alpha \subseteq 2^N | \forall X,Y \in \alpha : X \not\subseteq Y \ and\  Y \not\subseteq X\}
$.
We will denote the antichains by Greek minuscules. in case $N=\{1,2,...,n\}$, we will use $\mathcal{A}_n \equiv \mathcal{A}_N$.
$|N|$ is the {\it dimension} of $\mathcal{A}_N$. 
Antichains  $\alpha, \beta \in \mathcal{A}_N$ can be partly ordered as
$ \alpha \le \beta \Leftrightarrow \forall A \in \alpha : \exists B \in \beta : A \subseteq B$.
Supremum and infimum for this partial order are always defined and generated by
$
\alpha \vee \beta =  max(\alpha \cup \beta)
$, $
\alpha \wedge \beta =  max(\{A \cap B | A \in \alpha, B \in \beta\})
$,
where the max-operator on an arbitrary set of sets produces an antichain containing only the non dominated sets.
The operators $\wedge$ and $\vee$ are called meet and join respectively. 
$(\mathcal{A}_N,\wedge,\vee)$ is a completely distributive lattice.
For $\alpha, \beta \in \mathcal{A}_N$, intervals in this lattice are defined by
$[\alpha,\beta] = \{\chi \in \mathcal{A}_N| \alpha \le \chi \le \beta\}$.
It is clear that $[\alpha,\beta] \not= \emptyset$ if and only if $\alpha \le \beta$.
Moreover, if $\chi_1, \chi_2 \in [\alpha,\beta]$, then $\chi_1 \wedge \chi_2 \in [\alpha,\beta] $ and $\chi_1 \vee \chi_2 \in [\alpha,\beta]$. 
For $\alpha \le \beta$, the mapping $h_{[\alpha,\beta]} : \mathcal{A}_N \rightarrow [\alpha,\beta] : \chi \rightarrow \alpha  \vee (\chi \wedge \beta )$ is a lattice homomorphism. 
We denote the empty antichain $\{\}$ by  $\bot$ and the largest antichain $\{N\}$ by $\top$.
For $X \subseteq N$, the antichain of immediate subsets of X is denoted by $X^-$, so $X^- = \{X \backslash \{x\} | x \in X\}$.
For an antichain $\alpha$, we will denote $\alpha^- = \bigvee_{X \in \alpha} X^-$ and $\alpha^+ = \bigvee_{X \in 2^N, X^- \le \alpha} \{X\}$. 
We find $\emptyset^- =  \{\emptyset\}^- = \bot, \bot^+ = \{\emptyset\}$.\\

\section{Interval posets}
\label{sec:intervalposets}
For each interval $[\alpha,\beta]$,  the poset $\mathcal{P}_{[\alpha,\beta]} = (\{X \subseteq N | \alpha \vee \{X\} \in ]\alpha,\beta]\},\subset)$ spans the interval in the sense that each element in 
$[\alpha,\beta]$ is the join of $\alpha$ and an antichain of elements in $\mathcal{P}_{[\alpha,\beta]}$. 
Clearly $\mathcal{P}_{[\bot,\top]} = (2^N,\subset)$ spans $\mathcal{A}_N$.
We will call $\mathcal{P}_{[\alpha,\beta]}$ the underlying interval poset of $[\alpha,\beta]$.
This section presents four lemmata on the structure of interval posets.
Lemma \ref{lem:posetcharacteristic} characterizes interval posets and shows how to compute the bottom and the top of a spanned interval.
%
%
\begin{lemma}
\label{lem:posetcharacteristic}
For $S \subseteq 2^N$, $(S,\subset )$ is the underlying poset of an interval in $\mathcal{A}_N$ if and only if 
\[
\forall A_1,A_2 \in S:\forall X  \in 2^N: A_1 \subseteq X  \subseteq A_2 \Rightarrow X \in S
\]
In this case, the spanned interval is given by $[\bigvee_{X \in S} (X^- \backslash S),\bigvee_{X  \in S} \{X\}]$.
\proof
As follows immediately from the definition, every underlying poset of an interval satisfies this property.
Given a set $S \subseteq 2^N$ satisfying the property, it is easy to see that $S = \mathcal{P}_{[\bigvee_{X \in S} (X^- \backslash S),\bigvee_{X  \in S} \{X\}]}$.
$\boxempty$
\endproof
\end{lemma}
We will not distinguish in notation between $S$ and $(S,\subset)$.
The interval spanned by an interval poset $S$ is denoted by $\mathcal{I}_S$, thus $\mathcal{I}_{\mathcal{P}_{[\alpha,\beta]}} = [\alpha,\beta]$.
The structure of an interval poset is directly related to the structure of the spanned interval.
Isomorphic interval posets span isomorphic intervals.
The join of two intervals $[\alpha_1,\beta_1] \vee [\alpha_2,\beta_2]$ is the set of all antichains that can be written as
the join of two elements from each of the intervals, i.e. $[\alpha_1,\beta_1] \vee [\alpha_2,\beta_2] \equiv \{\chi_1 \vee \chi_2 | \chi_1 \in [\alpha_1,\beta_1],\ \chi_2 \in [\alpha_2,\beta_2] \}$.  It is easy to see that $[\alpha_1,\beta_1] \vee [\alpha_2,\beta_2] = [\alpha_1 \vee \alpha_2, \beta_1 \vee \beta_2]$. 
In case the decomposition is unique for each antichain in $[\alpha_1,\beta_1] \vee [\alpha_2,\beta_2]$, the join is called {\it direct} and is denoted as $[\alpha_1,\beta_1] \ovee  [\alpha_2,\beta_2]$.
The following three lemmata follow easily.
%
%
Lemma \ref{lem:posetdecomposition} shows how Interval posets can be decomposed if the respective constituents are incomparable:
\begin{lemma}
\label{lem:posetdecomposition}
Given two interval posets $S_1$ and $S_2$, such that no two sets $X_1 \in S_1$ and $X_2 \in S_2$ are in a subset/superset relationship, 
then the poset $S_1 \cup S_2$ is an interval poset, and we have $\mathcal{I}_{S_1 \cup S_2} = \mathcal{I}_{S_1} \ovee \mathcal{I}_{S_2}$. 
\end{lemma}
Lemma \ref{lem:posetredundancy} and \ref{lem:posetindistinguishability} show how to remove specific elements from the universe without changing the structure of the interval poset.
%
%
\begin{lemma}
\label{lem:posetredundancy}
Given a poset $S$ such that every element in $S$ contains a given set $A$ as a subset. Then the poset $\{X \backslash A|X \in S\}$ is isomorphic to $S$.
\end{lemma}
%
%
\begin{lemma}
\label{lem:posetindistinguishability}
Given a poset $S$ such that every element $X \in S$ that contains at least one element of a given set $A$ contains all elements of $A$. 
For $a \in A$, the poset $S' = \{(X \backslash A) \cup \{a\} |X \in S, A \subseteq X\} \cup \{X|X \in S, A \not\subseteq X\}$ is isomorphic to $S$.
\end{lemma}
\section{Intervals}
\label{sec:intervals}
The above lemmata have consequences for the structure of the spanned intervals.
The current section presents three lemmata that allow to reduce or decompose a given interval.
To ease notation, we will use two additional operators on antichains.
The first operator is the {\it span} of an antichain which is the union of its elements. For an antichain $\alpha$, we denote $\cup \alpha = \bigcup_{X \in \alpha} X$.
The second is the {\it direct product} of two antichains. For two antichains $\alpha, \beta$ with $(\cup \alpha) \cap (\cup \beta) = \emptyset$ the direct product is given by $\alpha \otimes\beta = \{A \cup B| A \in \alpha, B \in \beta\}$.
Lemma \ref{lem:intervalboundaryrendundancy} shows how common sets in top and bottom have no impact on the structure of the interval:
%
%
\begin{lemma}
\label{lem:intervalboundaryrendundancy}
For any two antichains $\alpha, \beta$, we have $\mathcal{P}_{[\alpha,\beta]} = \mathcal{P}_{[\alpha \backslash \beta, \beta \backslash \alpha]}$.
\end{lemma}
%
%
Lemma \ref{lem:intervalredundancy} builds on lemma \ref{lem:posetredundancy} to reduce top and bottom of an interval if a specific structure is recognized:
\begin{lemma}
\label{lem:intervalredundancy}
Given a set $A \subset N$ and two antichains $\chi' \le \chi \in \mathcal{A}_{N \backslash A}$.
Let $\beta =\{A\} \otimes \chi$
and $\alpha = \{A\} \otimes \chi' \vee A^- \otimes \chi$. The intervals $[\alpha,\beta]$ and $[\chi',\chi]$ are isomorphic.  
\proof
Any element  $X \in P_{[\alpha,\beta]}$ must be a subset of one of the elements of $\beta$, so let 
$X \subseteq A \cup B$ with $B \in \chi$. 
If $A \not\subseteq X$, then $X \subseteq A' \cup B$ for some $A' \in A^-$ and 
$\alpha \vee \{X\} = \alpha$, so that $X \not\in P_{[\alpha,\beta]}$.
So, $P_{[\alpha,\beta]}$ satifies the condition of lemma \ref{lem:posetredundancy}, and it is isomorphic to the poset $\{X \backslash A | X \in P_{[\alpha,\beta]}\} $.
We now find $P_{[\alpha,\beta]} = \{A \cup C | \exists B \in \chi: C \subseteq B \ and\ \not\exists B' \in \chi': C \subseteq B'\}$
and $\{X \backslash A | X \in P_{[\alpha,\beta]}\} = \{C | \exists B \in \chi: C \subseteq B \ and\ \not\exists B' \in \chi': C \subseteq B'\} = P_{[\chi',\chi]}$.
 $\boxempty$
\endproof
\end{lemma}
In Lemma \ref{lem:graphdecomposition}, we use the connected components of a graph associated with an interval to find a decomposition of this interval.
We associate a graph $\mathcal{G}_{[\alpha,\beta]}$ with each interval $[\alpha,\beta]$. The vertices of $\mathcal{G}_{[\alpha,\beta]}$ are the elements of  $\beta$.
Two vertices $A, B \in \beta$ are connected by an edge in  $\mathcal{G}_{[\alpha,\beta]}$ if and only if $\{A \cap B\} \not\le \alpha$. The connected components of $\mathcal{G}_{[\alpha,\beta]}$ determine the decomposability of the interval:
%
%
\begin{lemma}
\label{lem:graphdecomposition}
Let $\beta = \nu_1 \cup \nu_2$ with $\nu_1$ and $\nu_2$ such that no edge in the graph $\mathcal{G}_{[\alpha,\beta]}$ connects an element of $\nu_1$ with an element of $\nu_2$.
Then we have $[\alpha,\beta] = [\alpha, \alpha \vee \nu_1] \ovee [\alpha, \alpha \vee \nu_2] = [\alpha \wedge \nu_1,\nu_1] \ovee [\alpha \wedge \nu_2,\nu_2]$
\proof
Let $X \in \mathcal{P}_{[\alpha,\beta]}$, with $X \subseteq B_1 \in \nu_1$. Then $\not\exists B_2 \in \nu_2 : X \subseteq B_2$.
Indeed, if there were such a $B_2$, we would have $B_1 \cap B_2 \supseteq X$ and  $\{B_1 \cap B_2\} \not\le \alpha$ so that an edge in $\mathcal{G}_{[\alpha,\beta]}$
would connect $B_1 \in \nu_1$ with $B_2 \in \nu_2$.
The posets $\mathcal{P}_{[\alpha, \alpha \vee \nu_1]}$ and $\mathcal{P}_{[\alpha, \alpha \vee \nu_2]}$ are disconnected and the first part of the theorem follows.
To see the second part, note that $\alpha \wedge \nu_1  < \alpha \wedge \nu_1 \vee \{X\} \le \nu_1$ so that $X \in \mathcal{P}_{[\alpha \wedge \nu_1, \nu_1]}$. $\boxempty$
\endproof
\end{lemma}

\section{Partitioning by intervals}
\label{sec:partitioningbyintervals}
In this section, two theorems describe partitions of $\mathcal{A}_N$ or of intervals in $\mathcal{A}_N$ as disjoint unions of smaller intervals.
\subsection{Partitions associated with antichains}
\label{sub:antichainpartition}
In Theorem \ref{the:largestnondominatingdecomposition} (and Corollary \ref{cor:largestnondominatingdecomposition}) a partition is associated with each antichain in (an interval of)  $\mathcal{A}_N$.
To ease notation, we introduce the symbol $\check\chi$ for the largest antichain that does not dominate any of the sets in an antichain $\chi \in \mathcal{A}_N$:
\begin{eqnarray}
\label{eq:largestnondominating}
\check\chi & \in & \mathcal{A}_N \nonumber \\
\forall X \in \chi & : & \{X\} \not\le \check{\chi} \nonumber \\
\forall \rho \in \mathcal{A}_N : (\forall X \in \chi : \{X\} \not\le \rho) & \Rightarrow & \rho \le \check{\chi}
\end{eqnarray}
This maximum exists since the criterion is closed under the join operation and computing the largest non dominating antichain $\widecheck{(.)}$ is straightforward:
\begin{eqnarray*}
\forall A \subseteq N : \widecheck{\{A\}} & = & \{N \backslash \{a\} | a \in a\}\\
\forall \alpha,\beta \in  \mathcal{A}_N: \widecheck{\alpha \vee \beta} & = & \widecheck{\alpha} \wedge \widecheck{\beta}
\end{eqnarray*}
We can now prove the following
\begin{theorem}
\label{the:largestnondominatingdecomposition}
Let $\alpha \in \mathcal{A}_N$. We have
\[
 \mathcal{A}_N = \bigcup_{\chi \subseteq \alpha} [\chi,\widecheck{\alpha \backslash \chi}]
\]
and all intervals in the union are disjoint.
\proof
Let $\sigma \in  \mathcal{A}_N$ and let $\chi = (\sigma \wedge \alpha) \cap \alpha$.
Then $\sigma \ge \chi$ and $\sigma$ does not dominate any set in $\alpha \backslash \chi$ and hence $\sigma \in [\chi,\widecheck{\alpha \backslash \chi}]$.
Conversely, for $\sigma \in [\chi,\widecheck{\alpha \backslash \chi}]$, and $X \in  \sigma \wedge \alpha$ we have $X \in \alpha$. Since $\alpha \le \widecheck{\alpha \backslash \chi}$, $\{X\} \le  \widecheck{\alpha \backslash \chi}$
and $X \in \chi$. $\boxempty$
\endproof
\end{theorem}
\begin{corollary}
\label{cor:largestnondominatingdecomposition}
For $\alpha' \in [\alpha,\beta]$, we have
\[
 [\alpha,\beta] = \bigcup_{\chi \subseteq \alpha'} [\alpha \vee \chi,\beta \wedge \widecheck{\alpha' \backslash \chi}]
\]
and all intervals in the union are disjoint.
\begin{proof}
Note that, for general intervals, $[\alpha_1,\beta_1] \cap [\alpha_2,\beta_2] = [\alpha_1 \vee \alpha_2,\beta_1 \wedge \beta_2] $. $\boxempty$
\end{proof}
\end{corollary}

\subsection{Decomposition using the direct product}
\label{sub:compositiondirectproduct}
In \cite{PDCSDW2013, PDCSDW2011}, the direct product operator was used to decompose a space of a certain dimension according to spaces of lower dimension.
Here we derive this general decomposition from a lemma on interval posets.
\begin{lemma}
\label{lem:projectionlemma}.
Let $N = \{1,2,\ldots,n\}$ be a set of natural numbers and $\{N_1, N_2\}$ a partition of $N$.
Given two non-empty antichains $\alpha_1 \in \mathcal{A}_{N_1}$ and  $\alpha_2 \in \mathcal{A}_{N_2}$, let $\mathcal{S}_{\alpha_1,\alpha_2} = \{A_1 \cup A_2 |  \emptyset \subsetneq A_1 \subseteq N_1, \emptyset \subsetneq  A_2 \subseteq N_2, \{A_1\} \le \alpha_1,\{A_2\} \le \alpha_2\}$. $\mathcal{S}_{\alpha_1,\alpha_2} $ is an interval poset underlying the interval $[\alpha_1 \vee \alpha_2, \alpha_1 \otimes \alpha_2]$
\begin{proof}
Prove that $\mathcal{S}_{\alpha_1,\alpha_2} $ satisfies the conditions of lemma \ref{lem:posetcharacteristic} and note that 
\begin{itemize}
\item $\bigvee_{X \in \mathcal{S}_{\alpha_1,\alpha_2}} \{X\} = \alpha_1 \otimes \alpha_2$
\item $\bigvee_{X \in \mathcal{S}_{\alpha_1,\alpha_2}} \{X^{-} \backslash \mathcal{S}_{\alpha_1,\alpha_2}\} = \alpha_1 \vee \alpha_2$
\end{itemize}
$\boxempty$
\end{proof}
\end{lemma}
We use this lemma to prove the following theorem
\begin{theorem}
\label{the:intervalintervaldecomposition}
Let $N = \{1,2,\ldots,n\}$ be a set of natural numbers and $\{N_1, N_2\}$ a partition of $N$.
The set of intervals  
\[
\{[\bot,\bot]\} \cup \{[\alpha_1 \vee \alpha_2, \alpha_1 \otimes \alpha_2] | \alpha_1 \in \mathcal{A}_{N_1}, \alpha_2 \in \mathcal{A}_{N_2}, \alpha_1 \not= \bot, \alpha_2 \not= \bot\}
\] 
is a partition of $\mathcal{A}_N$.
\begin{proof}
Any non-empty antichain $\chi \in \mathcal{A}_{N}$ has a unique decomposition of the form $\alpha_1 \vee \alpha_2 \vee \alpha$ where $\alpha_1 \in \mathcal{A}_{N_1}, \alpha_2 \in \mathcal{A}_{N_2}$ and $\alpha$ is an anti-chain in $\mathcal{S}_{\alpha_1,\alpha_2}$ as defined in lemma \ref{lem:projectionlemma} satisfying 
$\alpha \wedge \{N_1\} \le \alpha_1, \alpha \wedge \{N_2\} \le \alpha_2$. Indeed, we find $\alpha_1 = \chi \wedge \{N_1\}, \alpha_2 = \chi \wedge \{N_2\}, \alpha = \chi \backslash (\alpha_1 \cup \alpha_2)$
and $\alpha_1 \vee \alpha_2 \le \alpha_1 \vee \alpha_2 \vee \alpha \le \alpha_1 \otimes \alpha_2$.
 $\boxempty$
\end{proof}
\end{theorem}
\section{Decomposition based on properies of antichains}
\label{sec:decomposebyproperty}
In this section, properties of antichains are used to derive decompositions or expressions for the size of (intervals in) $\mathcal{A}_N$.
An efficient expansion for the size of an interval in powers of two is given in Section \ref{sub:countingbysetsize}.
New numbers counting the antichains distinguishing between each of the elements of $N$ are introduced in Subsection \ref{sub:dedekindstirling2}.
These numbers are related to the Dedekind numbers through the  Stirling numbers of the second kind.
Subsection \ref{sub:connectedness} introduces connected antichains and relates their numbers to Dedekind numbers.
The numbers introduced in Subsections \ref{sub:dedekindstirling2} and \ref{sub:connectedness} can be computed from the known Dedekind numbers and are tabulated in Subsection \ref{sub:numericalvalues}.
Finally, an expansion of Dedekind numbers in products of interval sizes in a lower dimensional space is presented in \ref{sub:pcoefficients}.
The coefficients in this expansion are defined for any dimensions.
In the special case where the difference between the higher and lower dimensions is equal to two, the coefficients are powers of two and can be efficiently computed. 
\subsection{Summing over chains of constant set size}
\label{sub:countingbysetsize}
Antichains  of special interest in an interval are those existing of equally sized sets in the interval poset. 
Given an interval $[\alpha,\beta]$, with $\alpha < \beta$ so that $\mathcal{P}_{[\alpha,\beta]} \not= \emptyset$, we denote the size of the smallest sets in $\mathcal{P}_{[\alpha,\beta]}$ by $l_0$ and for the higher sizes $l_i = l_0 + i$.
For any $l \in \mathbb{N}$, define $\lambda_l([\alpha,\beta]) =   \{X \in \mathcal{P}_{[\alpha,\beta]} | |X| = l\}$.
In case there is no ambiguity, we will leave out the explicit reference to the interval $(\lambda_l([\alpha,\beta]) \equiv \lambda_l)$.
We will use $(.)^+$ and $(.)^-$ operators defined with respect to the interval $[\alpha,\beta]$.
So, for $\chi \subseteq \lambda_l$ we have
$\chi^+ = \{X \in \lambda_{l+1} | X^- \cap \lambda_l \subseteq \chi\}$, and for $l > l_0$,
$\chi^- = \{X^- \cap \lambda_{l-1} | X \in \chi\}$.
A crucial observation in this section is
\begin{lemma}
\label{lem:leveldecomposition}
Any antichain $\chi$ in an interval $[\alpha,\beta]$ has a unique decomposition as
\[
\chi = \alpha \vee \chi_0 \vee \chi_1 \vee \chi_2 \ldots
\]
where  $\forall i \in \mathbb{N}$, $\chi_i \in \lambda_{l_i}$ and $\chi_{i+1} \subseteq \chi_i^+$.
\proof
For each $i$, $\chi_i$ must at least contain the sets in $\chi$ of size $l_i$.
The constraint $\chi_{i+1} \subseteq \chi_i^+$ requires all subsets of size $i$ in $\mathcal{P}_{[\alpha,\beta]}$ of larger sets in $\chi$ to be present in $\chi_i$.
This implies $\chi_0 = (\chi \backslash \alpha) \wedge \lambda_{l_0}$ and $\forall i \in \mathbb{N}: \chi_{i+1} = (\chi \backslash (\alpha \vee \chi_0 \vee \chi_1 \ldots \vee \chi_i)) \wedge \lambda_{l_{i+1}}$.
$\boxempty$
\endproof
\end{lemma}

\begin{theorem}
\label{the:intervallatticesize}
For $\alpha \le \beta \in \mathcal{A}_N$, the size of the interval $[\alpha,\beta]$ is given by
\begin{eqnarray}
|[\alpha,\beta]| & = & \sum_{\chi_0 \subseteq \lambda_{l_0}} \sum_{\chi_2 \subseteq \chi_0^{++}}  \sum_{\chi_4 \subseteq \chi_2^{++}} \ldots 2^{|\chi_0^{+}| - |\chi_2^-| + |\chi_2^+| - |\chi_4^-| \ldots} 
\label{eqn:intervallatticesizeeven}\\
                          & = & \sum_{\chi_1 \subseteq \lambda_{l_1}} \sum_{\chi_3 \subseteq \chi_1^{++}}  \sum_{\chi_5 \subseteq \chi_3^{++}} \ldots 2^{|\lambda_{l_0}| - |\chi_1^{-}| + |\chi_1^{+}| - |\chi_3^-| + |\chi_3^+| - |\chi_5^-| \ldots}
\label{eqn:intervallatticesizeodd}
\end{eqnarray}
\proof
Due to lemma \ref{lem:leveldecomposition}, every tupple $(\chi_0,\chi_1,\ldots)$ satisfying the constraints in the lemma, corresponds with exactly one antichain in the interval.
In the current theorem, we sum over the possible assignments at even, respectively odd, levels.
Due to the lemma, $\chi_{i+2} \subseteq \chi_{i}^{++}$. The degree of freedom that is left after selection of $\chi_i$ and $\chi_{i+2}$ is given by the available sets at level $l_{i+1}$.
The sets in $\chi_{i}^+$ can be selected while the sets in $\chi_{i+2}^-$ must be selected. These degrees of freedom determine the power of 2 in the sum.
$\boxempty$
\endproof
\end{theorem}
Formulae \ref{eqn:intervallatticesizeeven} and \ref{eqn:intervallatticesizeodd} allow computing sizes of intervals up to $\mathcal{A}_6$ efficiently.
Theorem \ref{the:intervallatticesize}, as well as the decomposition Lemma \ref{lem:graphdecomposition} , were derived without relying on interval posets in \cite{PDCSDW2013}.
In the latter publication, the following example was given:
\begin{example}
\label{ex:latticesize}
Consider intervals of the form $I_N = [\{\emptyset\}, \binom{N}{2}]$ 
where, for convenience, we use the notation $\binom{N}{k}$ to denote the set of subsets of size $k$ of a set $N = \{1,2,...,n\}$.
$I_N$ has only two nonempty levels ($\lambda_0 = \bot$, $\lambda_{1}$ is the set of all singletons of elements of $N$ and $\lambda_2 = \binom{N}{2}$).
Applying formula \ref{eqn:intervallatticesizeodd}, we find $|I_N|  =  \sum_{i=0}^{n} \binom{n}{i} 2^{\binom{i}{2}}$ which is the well known formula for the number of  labeled graphs with at most $n$ nodes 
(Sloane series A006896 \cite{SLOANE}). This identity becomes obvious when we apply the alternative expression in formula \ref{eqn:intervallatticesizeeven} to find : 
$|I_N| =  \sum_{graphs \ g\ on\ n\ vertices} 2^{n - |vertices\ in\ g|} = \sum_{i=0}^{n} |graphs\ covering\ \{1,2,...,i\}| \binom{n}{i} 2^{n-i}$.
\end{example}
 
\subsection{Relation with Stirling numbers of the second kind}
\label{sub:dedekindstirling2}
This subsection considers intervals of the form 
$
[\{\{1\},\{2\},\ldots,\{n\}\},\top]]
$
in $\mathcal{A}_N$ 
which we will call {\it basic} intervals of dimension $n$. We will denote these by $\mathcal{B}_n$.
Within a basic interval, one can ask for those antichains in which the elements of a subset $X \subseteq N$ occur and are never separated.
Examples  of such antichains for $X = \{1,2\}$ in the basic interval of dimension 4 are $\{\{1,2\},\{3\},\{4\}\}$, $\{\{1,2,3\},\{1,2,4\}\}$ and $\{\{1,2,3,4\}\}$.
The antichain $\{\{1,2\},\{1,3\},\{2,4\}\}$ would not satisfy the criterion.
This set is spanned by a poset satisfying the conditions of  lemma \ref{lem:posetindistinguishability}.
It is not an interval, but due to lemma \ref{lem:posetindistinguishability}, it is isomorphic with, $\mathcal{B}_3 = [\{\{1\},\{2\},\{3\}\},\{\{1,2,3\}\}]$.
More generally, for any antichain of disjoint sets $\mathcal{\pi}$, denote by $\mathcal{B}_{\mathcal{\pi}} = \{\alpha \in \mathcal{A}_{\cup \mathcal{\pi}} |\pi \le \alpha \ and\  \forall A \in \alpha : \forall X \in \mathcal{\pi} : X \subseteq A \ or \ X \cap A = \emptyset\}$. For example, let $\pi = \{\{1,2\},\{3,4\}\}$. Then $\mathcal{B}_{\pi} = \{\{\{1,2\},\{3,4\}\},\{\{1,2,3,4\}\}\}$.
One  finds
\begin{lemma}
\label{lem:basicintervalcongruences}
For any antichain of disjoint sets $\pi$:
\[
 (\mathcal{B}_{\pi},\le) \cong (\mathcal{B}_{|\pi|},\le)
 \]
\end{lemma}
Given an antichain  $\alpha$, the partition  $\pi$ of $\cup \alpha$ of largest size such that $\alpha \in \mathcal{B}_{\pi}$ can be determined.
The sets in this partition are those for which the elements are never separated in $\alpha$, hence the following
\begin{definition}
Given $k \le n \in \mathbb{N}$, the distinguishing set $\mathcal{D}_{n,k}$ is the set of antichains $\alpha \in \mathcal{B}_n$ for which the largest partition $\pi$ of $N$ such that $\alpha \in \mathcal{B}_\pi$ has size $k$.
\end{definition}
These sets form a partition of the basic interval: $\mathcal{B}_n = \bigcup_{k=0}^{n} \mathcal{D}_{n,k}$
\footnote{$\mathcal{D}_{0,0} = \{\bot,\{\emptyset\}\}$ and for $n > 0$, we have $\mathcal{D}_{n,0} = \emptyset $.}.
Clearly, $\mathcal{D}_{n,k}$ can be further partitioned with respect to the partitions of $N$ of size $k$. 
Every component in this partition is isomorphic to $\mathcal{D}_{k,k}$. If we denote $\mathcal{D}_{k} \equiv \mathcal{D}_{k,k}$, we arrive at the following
\begin{lemma}
\label{lem:basicintervaldecomposition}
For $k \le n \in \mathbb{N}$, 
\begin{eqnarray}
|\mathcal{B}_n| & = & \sum_{k=0}^{n} |\mathcal{D}_{n,k}| \\
                            & = & \sum_{k=0}^{n} {n \brace k} |\mathcal{D}_{k}|
\end{eqnarray}
where ${n \brace k}$ is the Stirling number of the second kind. 
\begin{proof}
The lemma follows from the preceding observations since ${n \brace k}$ is the number of partitions of size $k$ of a set of size $n$.
\end{proof}
\end{lemma}
Given the well known decomposition for the Dedekind numbers $|\mathcal{A}_n| = \sum_{k = 0}^n \binom{n}{k}|\mathcal{B}_k|$, we arrive at
\begin{theorem}
\label{the:dedekindstirling2}
For $n \in \mathbb{N}$, the number of antichains of subsets of the set $\{1,2,\ldots,n\}$ is given by
\begin{eqnarray}
|\mathcal{A}_n|  & = & \sum_{k=0}^{n} {n+1 \brace k+1}  |\mathcal{D}_{k}|     
\end{eqnarray}
\begin{proof}
Given the previous discussion, and through lemma \ref{lem:basicintervaldecomposition}, we find
\begin{eqnarray*}
|\mathcal{A}_n|  & = &  \sum_{k = 0}^n \binom{n}{k}|\mathcal{B}_k| \\
             & = & \sum_{l=0}^{n} \binom{n}{l}  \sum_{k=0}^{l} {l \brace k}|\mathcal{D}_{k}|\\
             & = & \sum_{k=0}^{n} \big{(}\sum_{l=k}^{n} \binom{n}{l}  {l \brace k} \big{)}  |\mathcal{D}_{k}| \\
             & = & \sum_{k=0}^{n} {n+1 \brace k+1}  |\mathcal{D}_{k}|            
\end{eqnarray*} $\boxempty$
\end{proof}
\end{theorem}
Theorem \ref{the:dedekindstirling2} can be understood if one adds an element to $N$ that is never used and leave out the set in any partition containing this element.
Terms such as ${n+1 \brace k+1} |\mathcal{D}_k|$ then account for the number of partitions of which exactly $k$ components are used, not separated and not joined.

\subsection{Connectedness in antichains}
\label{sub:connectedness}
For a partition $\pi$ of $N$, an antichain  $\beta \in \mathcal{B}_{\pi}$ is built of unions of sets in $\pi$. 
If $X \in \pi$ is such a set, elements in $X$ are never separated in $\beta$.
$\beta$ thus defines an equivalence relation between elements in $N$ of which the sets in $\pi$ are the equivalence classes.
In general, any equivalence relation defined by an antichain generates a partition of $N$. 
In this subsection, we consider connectedness as another example of an equivalence relation.
Given an antichain $\beta \in \mathcal{B}_N$, consider the graph with vertices in $N$, and in which there is an edge between any two elements that are together in at least one set of $\beta$.
The connected components of this graph form a partition of $N$.
For a finite set $K \in \mathbb{N}_0$, we denote by $\mathcal{C}_K$ the set of antichains in $\mathcal{B}_K$ for which all elements of $K$ are in one such connected component.
It is clear that, for $k \ge 0$, all $\mathcal{C}_K$ for which $|K| = k$ are isomorphic and $|\mathcal{C}_K|$ only depends on $k$. As usual, we will use $\mathcal{C}_k$ and $\mathcal{C}_K$ interchangably in case no confusion is possible.
For $N = K_1 \cup K_2$ a partition of $N$, the number of antichains for which  $K_1$ and $K_2$ are connected is given by $|\mathcal{C}_{K_1}| * |\mathcal{C}_{K_2}|$. 
For $n = k_1 + k_2 + \ldots + k_p$, the number of antichains consisting of $p>0$ connected components with {\it different} sizes given by $k_1 > k_2 > \ldots > k_p$ is given by
$\binom{n}{k_1,k_2,\ldots,k_p} * |\mathcal{C}_{k_1}|*|\mathcal{C}_{k_2}|*\ldots*|\mathcal{C}_{k_p}|$.
In case sizes can be equal, $k_1 \ge k_2 \ge \ldots \ge k_p$, the multinomial coefficient does not account for interchangeability of whole sets in the partition.
In case $n = n_1 * k_1 + n_2 * k_2 + \ldots + n_p*k_p$, the number of antichains for which there are $n_1$ components of size $k_1$, $n_2$ components of size $k_2$, and so on, is given by 
\[
\frac{1}{n_1!n_2!\ldots n_p!}\binom{n}{k_1(n_1),k_2(n_2),\ldots,k_p(n_p)} 
* |\mathcal{C}_{k_1}|^{n_1}*|\mathcal{C}_{k_2}|^{n_2}*\ldots*|\mathcal{C}_{k_p}|^{n_p}
\]
where the notation $k(c)$ in the multinomial coefficient indicates $c$-fold repetition of the number $k$, e.g. $\binom{5}{2(2),1(1)} \equiv \binom{5}{2,2,1}$.
Since $C_{\emptyset} = \{\bot, \{\emptyset\}\}$,  $|\mathcal{C}_0| = 2$, leading to
\begin{theorem}
\label{the:connectedness}
For $n \in \mathbb{N}_0$,
\begin{equation}
\label{eq:connectedness}
|\mathcal{B}_n| = \sum_{n = \sum_{i=1}^{i=p}n_ik_i} \frac{1}{n_1!n_2!\ldots n_p!}\binom{n}{k_1(n_1),k_2(n_2),\ldots,k_p(n_p)} 
* |\mathcal{C}_{k_1}|^{n_1}*|\mathcal{C}_{k_2}|^{n_2}*\ldots*|\mathcal{C}_{k_p}|^{n_p}
\end{equation}
where the sum is over all expansions of $n$ as a linear combination of $p>0$ different values $k_1 > k_2 > \ldots > k_p > 0$.
\end{theorem}
Isolating the term in Theorem \ref{the:connectedness} corresponding to $p=1, n_1 = 1, k_1 = n$, we find 
\begin{eqnarray}
|\mathcal{B}_n| & = & |\mathcal{C}_n| \\ 
& + & \sum_{\substack{n = \sum_{i=1}^{i=p}n_ik_i \\ all\ k_i < n}} \frac{1}{n_1!n_2!\ldots n_p!}\binom{n}{k_1(n_1),k_2(n_2),\ldots,k_p(n_p)} 
* |\mathcal{C}_{k_1}|^{n_1}*|\mathcal{C}_{k_2}|^{n_2}*\ldots*|\mathcal{C}_{k_p}|^{n_p} \nonumber
\end{eqnarray}
relating $|\mathcal{C}_n|$ recursively to $|\mathcal{B}_n|$.
\begin{remark}
A similar recursion relation can be derived for any subset of $\mathcal{B}_n$ that is symmetric under permutation. 
Impose e.g. the restriction that the sets in the antichain be of some specific size $k$.
In this case, one replaces $|\mathcal{B}_n|$ by the total number of coverings of a set of size $n$ by subsets of size $k$.
For $k=2$, one thus finds a recursion relation for the number of connected graphs on $n$ nodes, as given in OEIS series A001187 (\cite{SLOANE}).
The results have been checked numerically \cite{CODESREPORTS}.
\end{remark}
\subsection{Numerical values}
\label{sub:numericalvalues}
Computing $|\mathcal{A}_n|$, $|\mathcal{B}_n|$, $|\mathcal{C}_n|$ and $|\mathcal{D}_n|$ for specific values of $n$ are all, given the analysis above, equivalent and equally hard. 
In the next section, we present powerful formulae for  $|\mathcal{A}_n|$ and $|\mathcal{B}_n|$. Presently, little is known about computing $|\mathcal{C}_n|$ and $|\mathcal{D}_n|$ directly.  Here we use
\begin{eqnarray}
|\mathcal{B}_n| & = & |\mathcal{A}_n| - \sum_{k=0}^{n-1} {n \choose k} |\mathcal{B}_k| \\
|\mathcal{D}_n| & = & |\mathcal{B}_n| - \sum_{k=1}^{n-1} { n \brace k} |\mathcal{D}_k| \\
|\mathcal{C}_n| & = & |\mathcal{B}_n| \\ 
& - & \sum_{\substack{n = \sum_{i=1}^{i=p}n_ik_i \\ all\ k_i < n}} \frac{1}{n_1!\ldots n_p!}\binom{n}{k_1(c_1),\ldots,k_p(c_p)} 
* |\mathcal{C}_{k_1}|^{n_1}*\ldots*|\mathcal{C}_{k_p}|^{n_p} \nonumber
\end{eqnarray}
to find the results in Table \ref{tbl:numbers}. 
The value of $|\mathcal{B}_n|$ could be computed independently using Theorem \ref{the:intervallatticesize} and the results for $|\mathcal{C}_n|$ and $|\mathcal{D}_n|$ have been checked by explicit enumeration up to $n = 6$.
\begin{landscape}
\begin{table}[htdp]
\caption{Dedekind numbers and related quantities}
\begin{center}
 \label{tbl:numbers}
\begin{tabular}{|c|r|r|r|r|}
\hline
$n$ & $\mathcal{A}_n$ & $\mathcal{B}_n$ & $\mathcal{C}_n$ & $\mathcal{D}_n $\\
\hline
\hline
0 & 2 & 2 & 2 & 2 \\
\hline
1 & 3 & 1 & 1 & 1 \\
\hline
2 & 6 & 2 & 1 & 1 \\
\hline
3 & 20 & 9 & 5 & 5 \\
\hline
4 & 168 & 114 & 84 & 76 \\
\hline
5 & 7581 & 6894 & 6348 & 5993 \\
\hline
6 & 7828354 & 7785062 & 7743728 & 7689745 \\
\hline
7 & 2414682040998 & 2414627396434 & 2414572893530 & 2414465044600 \\
\hline
8 & 56130437228687557907788 & 56130437209370320359966 & 56130437190053299918162 & 56130437141763247212112 \\
\hline
\end{tabular}
\end{center}
\label{default}
\end{table}
\end{landscape}
%

\subsection{P-coefficients}
\label{sub:pcoefficients}
For $n \in \mathbb{N}, N = \{1,\ldots,n\}, N_1 \subsetneq N$, clearly the corresponding posets satisfy $\mathcal{P}_{[\bot,\{N_1\}]} \subsetneq \mathcal{P}_{[\bot,\top]}$
and any set in  $\mathcal{P}_{[\bot,\top]}$ is a superset of at least one set in $ \mathcal{P}_{[\bot,\{N_1\}]}$.
Note that $\mathcal{P}_{[\bot,\top]} = 2^N$.
As a consequence, we have
\begin{lemma}
\label{lem:uniquedecomposition}
Let $n \in \mathbb{N}, N = \{1,\ldots,n\}, N_1 \subsetneq N$.
Each $\chi \in \mathcal{A}_n$ has a unique decomposition as
\begin{equation}
 \label{lem:decompositionformula} 
\chi = \bigvee_{P \in  2^{N_1}} \chi_P \times \{P\}
\end{equation}
where for $P, P' \in 2^{N_1}$, $\chi_P,\chi_{P'} \in \mathcal{A}_{N \backslash N_1}$ and
\begin{equation}
\label{lem:uniqueness condition}
P \subseteq P' \Rightarrow \chi_P \ge \chi_{P'}
\end{equation}
\begin{proof}
For each $P \subseteq N_1$, let  $\chi_P = \{X \backslash N_1 | X \in \chi, P \subseteq X\}$.
It is easy to see that the resulting  tuple $(\chi_P | P \subseteq N_1)$ satisfies the conditions of the lemma.
To see the uniqueness, derive from equation \ref{lem:decompositionformula} that $\chi_{N_1} =  \{X \backslash N_1 | X \in \chi, N_1 \subseteq X\}$.
Remove $\chi_{N_1} \times \{N_1\}$ from $\chi$ and derive that $\chi_P = \{X \backslash N_1 | X \in \chi, P \subseteq X\}$ for each immediate subset $P \subset N_1$ 
from equation \ref{lem:decompositionformula}  using condition \ref{lem:uniqueness condition}. Repeating recursively yields the result. $\boxempty$
\end{proof}
\end{lemma}
$\mathcal{A}_n$ is hence isomorphic to the set of order reversing mappings $2^{N_1} \rightarrow \mathcal{A}_{(n-|N_1|)}$.
Using this notation, given $\rho_0 \ge \rho_{N_1}$, one can now ask for the number of antichains $\chi$
for which $\bigwedge_{i \in N_1} \chi_{\{i\}} = \rho_{N_1}$ and $\bigvee_{i \in N_1} \chi_{N_1 \backslash \{i\}} = \rho_0$.
We will address this number as the $\mathcal{P}-coefficient$ :
\begin{definition}
\label{def:pcoefficients}
For $n,k \in \mathbb{N}$, $\rho_1, \rho_2 \in \mathcal{A}_n$, let $K = \{n+1,\ldots,n+k\}$. 
$\mathcal{P}_{n,k,\rho_1, \rho_2}$ is the number of solutions $(\chi_P \in \mathcal{A}_n | P \subseteq K)$ to the simultaneous equations
\begin{eqnarray}
\chi_{\emptyset} & = & \rho_2 \\
\chi_K & = & \rho_1 \\
\bigwedge_{i \in K} \chi_{K \backslash \{i\}} & = & \rho_1
 \label{def:pcoefflowerconstraint}\\
\bigvee_{i \in K}  \chi_{\{i\}}  & = & \rho_2 
\label{def:pcoeffupperconstraint} \\
\forall P,P' \in 2^K, P \subseteq P' & \Rightarrow & \chi_{P} \ge \chi_{P'}
\label{def:pcoefforderreversal}
\end{eqnarray}
\end{definition} 
One finds easily, $\mathcal{P}_{n,k,\rho_1, \rho_2} = 0$ for $\rho_1 > \rho_2 \in \mathcal{A}_n$,  $\mathcal{P}_{n,k,\rho, \rho} = 1$ for any $n, k \in \mathbb{N}$ and 
any $\rho \in \mathcal{A}_n$, $\mathcal{P}_{n,0,\rho_1, \rho_2} = 0$ for $\rho_1 \not= \rho_2$ ($\chi_K = \chi_{\emptyset}$), $\mathcal{P}_{n,1,\rho_1, \rho_2} = 0$ for $\rho_1 \not= \rho_2$ ($\bigwedge \chi_{\emptyset} = \bigwedge \rho_1 = \rho_2$). $\mathcal{P}_{n,k,\bot, \{\emptyset\}} = |\mathcal{A}_k| - 2$ for $k > 0$  (the equations are equivalent with the definition of monotone boolean functions on $2^K$ taking $false$ on $\emptyset$ and $true$ on $K$), An important non-trivial case is $\mathcal{P}_{n,2,\rho_1, \rho_2}$ for $\rho_1 \le \rho_2$. 
\begin{lemma}
Let $n \in \mathbb{N}, \rho_1 \le \rho_2 \in \mathcal{A}_n$. Then we have
\begin{equation}
\mathcal{P}_{n,2,\rho_1,\rho_2} = 2^{|C_{\rho_1,\rho_2}|}
\end{equation} 
\begin{proof}
\end{proof}
Let $A, B \in \rho_2$ such that $\{A \cap B\} \not\le \rho_1$. 
Equation \ref{def:pcoeffupperconstraint} implies that either $A \in \chi_{\{1\}}$ or $A \in \chi_{\{2\}}$ and equivalently for $B$.
Equation \ref{def:pcoefflowerconstraint} however, forbids $A \in \chi_{\{1\}}$ and $B \in \chi_{\{2\}}$ so $A$ and $B$ must be in exactly one of $\chi_{\{1\}}$ and $\chi_{\{2\}}$.
By extension, each connected component in $C_{\rho_1,\rho_2}$ must be in exactly one of $\chi_{\{1\}}$ and $\chi_{\{2\}}$. $\boxempty$
\end{lemma}
The importance of P-coefficients becomes apparent in the following
\begin{theorem}
\label{the:pcoefficientstheorem}
For $n, k \in \mathbb{N}, k > 1$
\begin{equation}
|\mathcal{A}_{n+k}| = \sum_{ \substack{\alpha,\beta \in \mathcal{A}_n\\\alpha \le \beta}}|[\bot,\alpha]|\mathcal{P}_{n,k,\alpha,\beta}|[\beta,\top]|
\end{equation}
\begin{proof}
Let $K = \{1,\ldots,n\}$. For $k > 1$, the tuple $(\chi_P | \emptyset \subsetneq P \subsetneq K)$ in the unique decomposition of lemma \ref{lem:uniquedecomposition} is not empty.
Given a  specific value of this tuple, $\chi_{K}$ can take any value from $[\bot,\bigwedge_{i \in K} \chi_{\{i\}}]$ and similarly $\chi_{\emptyset}$ can take values in 
$[\bigvee_{i \in K} \chi_{K \backslash \{i\}},\top\}]$. Given $\alpha \le \beta \in \mathcal{A}_n$, the number of such tuples for which $\bigwedge_{i \in K} \chi_{\{i\}} = \alpha$ and $\bigvee_{i \in K} \chi_{K \backslash \{i\}} = \beta$ is given by $\mathcal{P}_{n,k,\alpha,\beta}$. $\boxempty$
\end{proof}
\end{theorem}
The sum for $n=6$ still contains $|\mathcal{A}_6|^2 \approx 6,1\times10^{24}$ terms and involves $2*|\mathcal{A}_6| \approx 1,5\times10^7$ interval sizes.
Using permutation symmetry, the number of terms can be reduced to approximately $1,1\times10^{11}$ terms and $1,6\times10^4$ (the term for $n=6$ in series A003182 at \cite{SLOANE}) interval sizes.
Given the efficient algorithm for the $p$-coefficients of order $2$, the resulting formula allows computing $|\mathcal{A}_{8}$ in about 8 hours on a 4 core 3 Ghz computer. From the numbers it is clear that so far, $|\mathcal{A}_{9}$ is out of reach.

\section{Acknowledgement}
Work supported by the Belgian Science Policy Office (BELSPO) in the Interuniversity Attraction Pole COMEX. (http://comex.ulb.ac.be)

\section*{References}

\bibliography{mybibfile}

\begin{thebibliography}{10}
\expandafter\ifx\csname url\endcsname\relax
  \def\url#1{\texttt{#1}}\fi
\expandafter\ifx\csname urlprefix\endcsname\relax\def\urlprefix{URL }\fi
\expandafter\ifx\csname href\endcsname\relax
  \def\href#1#2{#2} \def\path#1{#1}\fi

\bibitem{DEDEKIND}
R.~Dedekind, Uber zerlegungen von zahlen durch ihre grossten gemeinsamen
  teiler, Gesammelte Werke 2 (1897) 103--148.

\bibitem{SLOANE}
N.~Sloane, The on-line encyclopedia of integer sequences. (oeis),
  http://www.research.att.com/~njas/sequences/.

\bibitem{KLEITMAN}
D.~Kleitman, On dedekind's problem: the number of isotone boolean functions,
  Proc. Amer. Math. Soc. 21 (1969) 677--682.

\bibitem{KLEITMAN_MARKOWSKY}
D.~Kleitman, G.~Markowsky, On dedekind's problem: the number of isotone boolean
  functions. ii, Trans. Amer. Math. Soc. 213 (1975) 373--390.

\bibitem{KORSHUNOV}
A.~D. Korshunov, The number of monotone boolean functions (russian), Problemy
  Kibernet. 38 (1981) 5--108.

\bibitem{KAHN}
J.~Kahn, Entropy, independent sets and antichains: a new approach to dedekind's
  problem, Proc. Amer. Math. Soc. 130 (2) (2002) 371--378.

\bibitem{CHURCH1940}
R.~Church, Numerical analysis of certain free distributive structures, Duke.
  Math. J. 6 (1940) 732Ð734.

\bibitem{CHURCH1965}
R.~Church, Enumeration by rank of the elements of the free distributive lattice
  with 7 generators, Notices Amer. Math. Soc. 12 (1965) 724.

\bibitem{WARD}
M.~Ward, Note on the order of free distributive lattices, Bull. Amer. Math.
  Soc. 52 (1946) 423.

\bibitem{WIEDEMANN}
D.~Wiedemann, A computation of the eighth dedekind number, Order 8 (1) (1991)
  5--6.

\bibitem{BAKOEV}
V.~Bakoev, One more way for counting monotone boolean functions, in: Thirteenth
  International Workshop on Algebraic and Combinatorial Coding Theory, 2012, p.
  47Ð52.

\bibitem{FYDITEK}
R.~Fidytek, A.~W. Mostowski, R.~Somla, A.~Szepietowski, Algorithms counting
  monotone boolean functions, Information Processing Letters 79 (2001) 203Ð209.

\bibitem{TAMON2014}
T.~Stephen, T.~Yusun,
  \href{http://dx.doi.org/10.1016/j.dam.2013.11.015}{Counting inequivalent
  monotone boolean functions}, Discrete Appl. Math. 167 (2014) 15--24.
\newblock \href {http://dx.doi.org/10.1016/j.dam.2013.11.015}
  {\path{doi:10.1016/j.dam.2013.11.015}}.
\newline\urlprefix\url{http://dx.doi.org/10.1016/j.dam.2013.11.015}

\bibitem{CODESREPORTS}
P.~De~Causmaecker, Codes reports, www.kuleuven-kortrijk.be/codes/.

\bibitem{PDCSDW2013}
P.~De~Causmaecker, S.~De~Wannemacker, Decomposition of intervals in the space
  of anti-monotonic functions, in: M.~Ojeda-Aciego, J.~Outrata (Eds.),
  Proceedings of the Tenth International Conference on Concept Lattices and
  Their Applications (CLA 2013), Vol. 1062, CLA 2013, La Rochelle, France,
  2013, pp. 57--67.

\bibitem{PDCSDW2011}
P.~De~Causmaecker, S.~De~Wannemacker, Partitioning in the space of anti
  monotonic functions, arXiv:1103.2877 [math.NT].

\end{thebibliography}

\end{document}